\documentclass[12pt,a4paper,reqno]{amsart}
\usepackage{amsmath}
\usepackage{amsfonts}
\usepackage{amssymb}
\usepackage{amscd}

\newcommand\rank{{\operatorname{rank}}}

\newcommand\R{{\mathbf{R}}}
\newcommand\C{{\mathbf{C}}}
\newcommand\Q{{\mathbf{Q}}}

\renewcommand\P{{\mathbf{P}}}
\newcommand\E{{\mathbf{E}}}
\newcommand\bv{{\mathbf{v}}}
\newcommand\bw{{\mathbf{w}}}

\newcommand\Z{{\mathbf{Z}}}

\newcommand\eps{\varepsilon}

\theoremstyle{plain}
  \newtheorem{theorem}[subsection]{Theorem}

  \newtheorem{proposition}[subsection]{Proposition}
  
  \newtheorem{lemma}[subsection]{Lemma}
  \newtheorem{corollary}[subsection]{Corollary}

\theoremstyle{remark}
  \newtheorem{remark}[subsection]{Remark}

\theoremstyle{definition}
  \newtheorem{definition}[subsection]{Definition}

\include{psfig}

\begin{document}

\title[Sharp Inverse Littlewood-Offord]{A sharp inverse Littlewood-Offord theorem}

\author{Terence Tao}
\address{Department of Mathematics, UCLA, Los Angeles CA 90095-1555}
\email{tao@math.ucla.edu}
\thanks{T. Tao is supported by a grant from the Macarthur Foundation.}

\author{Van Vu}
\address{Department of Mathematics, Rutgers, Piscataway, NJ 08854}
\email{vanvu@math.rutgers.edu-}
\thanks{V. Vu is supported by research grants DMS-0901216 and AFOSAR-FA-9550-09-1-0167.}


\begin{abstract}
Let $\eta_i, i=1,\dots, n$ be iid Bernoulli random variables. Given
a multiset $\bv$ of $n$ numbers $v_1, \dots, v_n$, the \emph{concentration
probability} $\P_1(\bv)$ of $\bv$ is defined as $\P_1(\bv) := \sup_{x} \P( v_1 \eta_1+ \dots
v_n \eta_n=x)$. A classical result of Littlewood-Offord and Erd\H os
from the 1940s asserts that if the $v_i $ are non-zero, then this
probability is at most $O(n^{-1/2})$.  Since then,  many researchers
obtained better bounds by assuming various restrictions on $\bv$.

In this paper, we  give an asymptotically optimal characterization
for all multisets  $\bv$ having  large concentration probability.
This allow us to strengthen or recover several previous results in a
straightforward manner.
\end{abstract}
\maketitle

\section{Introduction}

The purpose of this paper is to study the \emph{Littlewood-Offord} and \emph{inverse Littlewood-Offord} problems regarding concentration of random walks in torsion-free abelian groups.  We recall some notation from \cite{TVcond}.

\begin{definition}[Concentration probabilities]
Let $G = (G,+)$ be an additive group (e.g. the integers $\Z$, the complex numbers $\C$, or a vector space $\R^m$).
Let $\bv = (v_1,\ldots,v_n)$ be a multiset of $n$ elements of $G$ (allowing
repetitions). For any $0 \leq \mu \leq 1$, we define the \emph{lazy random walk} $S^{\mu}(\bv)$ with steps $\bv$ and density $\mu$ to be the $G$-valued random variable
$$S^{\mu}(\bv)  := v_1 \eta^{\mu}_1 + \dots + v_n \eta^{\mu}_n $$
\noindent where the
$\eta^{\mu}_i$'s  are iid copies of the (lazy coin flip) random
variable $\eta^{\mu}$
  which equals $0$ with probability $1-\mu$ and
$\pm 1$ with probability $\mu/2$ each.  We define the \emph{concentration probability}
$\P_\mu(\bv)$ to be the quantity
\begin{equation} \label{equa:defofPmu}  \P_\mu(\bv) :=
\max_{a \in G} \P( S^{\mu} (\bv) = a).\end{equation}
\end{definition}

\begin{remark} We are interested in the regime when $0 < \mu \leq 1$ is fixed and $n$ is large.  The most interesting case is perhaps when $\mu=1$. In this case $\eta$ is the Bernoulli random variable (fair coin flip), and $\P_1(\bv)$ is the maximum multiplicity among the $2^n$ signed sums $\pm v_1 \pm \ldots \pm v_n$, divided by $2^n$.
Such probabilities appear in many situations in combinatorics and the theory of random structures, for instance in understanding the singularity probability of discrete random matrices (see e.g. \cite{KKS}, \cite{TVdet}, \cite{TVsing}, \cite{RV}, \cite{TVcir}).
\end{remark}

We will assume throughout this paper that $G$ is \emph{torsion-free}, thus $nx \neq 0$ whenever $x \in G$ is non-zero and $n$ is a non-zero integer. In this case we can usually reduce to the model case $G=\Z$ by means of Freiman isomorphisms (see \cite[Lemma 5.25]{TVbook}).  

Broadly speaking, we expect $\P_\mu(\bv)$ to be large if and only if $\bv$ has significant additive structure.  To explore this phenomenon, we ask the following two general (and closely related) questions:

\begin{itemize}
\item (Forward Littlewood-Offord problem) Given additive structural hypotheses on $v_1,\ldots,v_n$, what bounds can one give for $\P_\mu(\bv)$?
\item (Inverse Littlewood-Offord problem) Given bounds on $\P_\mu(\bv)$, what can one say about the additive structure of the $v_1,\ldots,v_n$?
\end{itemize}

Let us now recall some previous results on these problems; further discussion may be found in \cite[Chapter 5]{TVbook}.  For simplicity we take $\mu=1$.  With no assumptions on $\bv = (v_1,\ldots,v_n)$, we easily obtain the inequalities
$$ 2^{-n} \leq \P_1(\bv) \leq 1$$
with the upper bound being attained precisely when all the $v_i$ are zero, and the lower bound attained precisely when the the $v_i$ are \emph{dissociated} (which means that all the $2^n$ partial sums $\sum_{i \in A} v_i$ with $A \subset \{1,\ldots,n\}$ are distinct).  These two cases represent extreme additive structure and extreme lack of additive structure respectively.

Throughout this paper we adopt the following asymptotic notation:

\begin{definition}[Asymptotic notation]  The asymptotic notation $X = O(Y)$, $X \ll Y$, $Y = \Omega(X)$, or $Y \gg X$ denotes the bound $X \leq CY$ for all $n \geq C$ and some absolute constant $C$; we also use $X = \Theta(Y)$ for $X \ll Y \ll X$.  Subscripting such as $O_d(Y)$ means that the implied constants $C$ in the asymptotic notation are allowed to depend on $d$.
\end{definition}

Littlewood and Offord \cite{LO} and then Erd\H{o}s \cite{Erd1} were able to improve the upper bound assuming that some of the $v_i$ were non-zero.  In particular, from the results in \cite{Erd1} one obtains the inequality
\begin{equation}\label{erdos-lo}
\P_1(\bv) \ll n^{-1/2}
\end{equation}
if all of the $v_i$ are non-zero (Littlewood and Offord obtained the slightly weaker bound $\P_1(\bv) \ll n^{-1/2} \log n$).   This bound is sharp: if $v_1=\ldots=v_n$, one easily verifies that $\P_1(\bv) \gg n^{-1/2}$ (and in fact this example gives the precise maximum value of $\P_1(\bv)$.

The above result is phrased as a forward Littlewood-Offord result, but can be easily rephrased as an inverse theorem:

\begin{theorem}[Erd\H{o}s' inverse Littlewood-Offord theorem]\label{eilo}  Let $\bv = (v_1,\ldots,v_n)$ be an $n$-tuple in a torsion-free additive group $G$.  Suppose that $\P_1(\bv) \gg k^{-1/2}$ for some $k \geq 1$.  Then all but $O(k)$ of the $v_1,\ldots,v_n$ are zero.
\end{theorem}

One can improve the upper bounds further by excluding the above counter-example.   Indeed, from the work of Erd\H{o}s and Moser \cite{erdmos} and then S\'ark\"ozy and Szemer\'edi \cite{SS}, the bound
\begin{equation}\label{ss-lo}
\P_1(\bv) \ll n^{-3/2}
\end{equation}
was established if all the $v_i$ were distinct (the earlier paper \cite{erdmos} establishes the slightly weaker bound $\P_1(\bv) \ll n^{-3/2} \log n$).
Again, this result is sharp, since if one takes $v_1,\ldots,v_n$ to be a proper arithmetic progression, one easily verifies that $\P_1(\bv) \gg n^{-3/2}$.  Later, Stanley \cite{Stan}, using algebraic methods, gave a very explicit bound for the optimal value of $\P_1(\bv)$.

The higher dimensional version of the problem, in which $G$ is a vector space $\R^m$, has also attracted attention.
Without the assumption that the $v_i$'s are different, the best
result was obtained by Frankl and F\"uredi in \cite{FF}, following earlier results by Katona \cite{Kat},
Kleitman \cite{Kle}, Griggs, Lagarias, Odlyzko and Shearer
\cite{GLOS} and many others. However, the techniques used in these
papers did not seem strong enough to recover \eqref{ss-lo}.  On the other hand, Hal\'asz \cite{Hal}, using harmonic analysis methods, managed to generalise \eqref{ss-lo}, proving even stronger bounds upon forbidding more additive correlations among the $v_i$'s. 

\begin{theorem}[Halasz inequality]\label{theorem:halasz}\cite{Hal}, \cite[Exercise 7.2.8]{TVbook} Let $\bv = (v_1,\ldots,v_n)$ be an $n$-tuple in a torsion-free additive group $G$.  Let $l \geq 1$ be an integer and let $0 < \mu \leq 1$.
Let $R_l$ be the number of solutions of the equation
$$\epsilon_1 v_{i_1} + \dots  +\epsilon_{2l} v_{i_{2l}} =0 $$
\noindent where $\epsilon_i \in \{-1,1\}$ and $i_1, \dots, i_{2l}$ are
(not necessarily different) elements of $\{1,2, \dots, n \}$. Then
$$\P_{\mu} (\bv) \ll_{l,\mu} n^{-2l-1/2} R_l. $$
\end{theorem}

It is easy to see that the $l=1$ case of Theorem \ref{theorem:halasz} implies the bound \eqref{ss-lo}.

\subsection{Main results}

Theorem \ref{theorem:halasz} states, roughly speaking, that if $\P_\mu(\bv)$ is large, then there is a large amount of additive structure (in the form of short additive relations) between the $v_i$.  Now we consider a slightly different type of additive structure, namely containment in a (symmetric) \emph{generalized arithmetic progression} (or GAP); we recall this concept in Section \ref{gap-appendix}.  It is not hard to show that if all the $v_i$ are contained in a GAP of bounded rank and controlled size, then the concentration probability $\P_\mu(\bv)$ is large.  More precisely, one has

\begin{proposition}[Forward Littlewood-Offord theorem]\label{filo}  Let $Q$ be a symmetric GAP in an additive group $G$ with rank $d$, and let $\bv = (v_1,\ldots,v_n)$ be such that $v_1,\ldots,v_n \in Q$.  Let $0 < \mu \leq 1$.  Then we have
$$ \P_\mu(\bv) \gg_d |Q_{\sqrt{\mu n}}|^{-1} \gg_d (1 + \mu n)^{-d/2} |Q|^{-1},$$
where the dilate $Q_t$ of $Q$ is defined in Section  \ref{gap-appendix}.
\end{proposition}

We prove this easy result in Section  \ref{gap-appendix}; it reflects the intuition that a lazy random walk with steps in $Q$ should mostly take values in the dilate $Q_{O(\sqrt{\mu n})}$.  Note that this result incorporates the examples used to demonstrate that \eqref{erdos-lo} and \eqref{ss-lo} are sharp.  See also \cite[Theorem 6.6]{TVcir} for a more complicated result in a similar spirit.

We now turn to the question of whether a converse to Proposition \ref{filo} exists.  In \cite{TVcond}, the authors showed

\begin{theorem}[Weak Inverse Theorem]\label{theorem:weak} Let $A, \eps > 0$ and $0 < \mu \leq 1$, and let $\bv = (v_1,\ldots,v_n)$ be an $n$-tuple in a torsion-free additive group $G$ be such that
$$\P_{\mu} (v) \ge n^{-A}. $$
Then there exists a proper symmetric GAP $Q$ of rank $d$ for some $d=O_{A,\eps}(1)$, of volume $O_{A,\mu,\eps}(n^B)$ for some $B = O_{A,\eps}(1)$, which contains all but $O_{A,\mu,\eps}(n^{1-\eps})$ elements of $\bv$ (counting multiplicity).
\end{theorem}

The reason we call Theorem \ref{theorem:weak} a \emph{weak} inverse theorem because the dependence of $B$ on $A$ is not optimal ($B$ is roughly $2A^2$).  The first main result of this paper is to obtain a sharper converse to Proposition \ref{filo}, in which $B$ is taken to be $A-\frac{d}{2}+\eps$:

\begin{theorem}[Strong Inverse Theorem] \label{theorem:strong}
Let $A, \eps > 0$, and let $\bv = (v_1,\ldots,v_n)$ be an $n$-tuple in a torsion-free additive group $G$ be such that
\begin{equation}\label{pmuv}
\P_{\mu} (v) \ge n^{-A}. 
\end{equation}
Then there exists a proper symmetric GAP $Q$ of rank $d \leq 2A$ of volume $O_{A, \mu, \eps}(n^{A - \frac{d}{2} + O_A(\eps)})$, which contains all but $O_{A,\mu,\eps}(n^{1-\eps} )$ elements of $\bv$ (counting multiplicity).
\end{theorem}

Comparing this with Proposition \ref{filo} we see that except for epsilons, the exponent $A-\frac{d}{2}+O_A(\eps)$ here cannot be improved.  

Theorem \ref{theorem:strong} will be deduced as the special case of the following stronger result.

\begin{theorem}[General Strong Inverse Theorem] \label{theorem:stronggeneral}
Let $d \geq 1$ be an integer  and let $0 < \eps, \mu < 1$ be
constants.  Then there is a constant $C_{0}=C_{0} (d, \eps, \mu)$ such
that the following holds for all sufficiently large $n$ and $k$
with $1 \leq k < \sqrt n$. Suppose that $\bv= (v_{1}, \dots, v_{n})$ is an $n$-tuple in a torsion-free additive group $G$ that satisfies
\begin{equation}\label{pi2}
\P_{\mu}(\bv) \ge C_{0}k^{-d}.
\end{equation}
Then there exists a proper symmetric GAP $Q$ of rank at most $d-1$ and volume
\begin{equation}\label{qde}
\operatorname{vol}(Q)  \le    \P_{\mu}(\bv)^{-1} k^{\eps}
\end{equation}
such that $Q_{1/k}$ contains all but at most $O_{d, \mu, \eps} (k^2
\log k)$ of the $v_1,\ldots,v_n$.  Furthermore, there is a positive integer $C=
C(d, \mu, \eps)$ such that the steps of $Q$ lie in $\{ v_1/C, \ldots, v_n/C \}$.
\end{theorem}

Let us see how this theorem implies Theorem \ref{theorem:strong}.

\begin{proof}[Proof of Theorem \ref{theorem:strong} assuming Theorem \ref{theorem:stronggeneral}]  Let $A, \mu, \eps, \bv, n, G$ be as in Theorem \ref{theorem:strong}.  By shrinking $\eps$ if necessary we may assume that $\eps$ is small depending on $A$.  We may assume that $n$ is large depending on $A,\mu,\eps$ as the claim is trivial otherwise. Let $d$ be the first integer larger than $2A$, and let $C_0$ be as in Theorem \ref{theorem:stronggeneral}.   For $\eps$ small and $n$ large, we see from \eqref{pmuv} that \eqref{pi2} holds for $k := n^{1/2 - \eps}$.  By Theorem \ref{theorem:stronggeneral}, we obtain a proper symmetric GAP $Q$ of rank $r$ at most $d-1$ and volume $O( n^{A+\eps} )$ such that $Q_{1/k}$ contains all but $O( n^{1-\eps} )$ of the $v_1,\ldots,v_n$.  Observe that any dimension of $Q$ that is less than $k$ does not contribute anything to $Q_{1/k}$, so by deleting these steps (and reducing the rank $r$ of $Q$) we may assume that all dimensions of $Q$ are at least as large as $k$.  Then $Q_{1/k}$ is a proper symmetric GAP of rank at most $2A$ and volume $O_{A,\mu,\eps}( k^{-r} |Q| ) = O_{A,\mu,\eps}( n^{A - r/2 + O_A(\eps)} )$, and the claim follows.
\end{proof}

\subsection{Applications}

We now give some applications of Theorem \ref{theorem:strong} and Theorem \ref{theorem:stronggeneral}.  We first observe that these theorems can recover the classical bounds \eqref{erdos-lo}, \eqref{ss-lo} except for epsilon losses:

\begin{proposition}\label{class}  Let $\bv = (v_1,\ldots,v_n)$ be an $n$-tuple in a torsion-free additive group $G$, and let $\eps > 0$ and $0 < \mu \leq 1$.
\begin{itemize}
\item[(i)] If all the $v_i$ are non-zero, then $\P_\mu(\bv) \ll_{\mu,\eps} n^{-1/2+\eps}$.
\item[(ii)] If all the $v_i$ are distinct, then $\P_\mu(\bv) \ll_{\mu,\eps} n^{-3/2+\eps}$.
\end{itemize}
\end{proposition}

\begin{proof} We may assume that $n$ is large compared to $\mu,\eps$, as the claim is trivial otherwise.

We first prove (i).  Suppose for contradiction that $\P_\mu(\bv) \geq n^{-1/2+\eps}$.  Applying Theorem \ref{theorem:strong} with $A := 1/2-\eps$ we see that there exists a symmetric GAP $Q$ of rank at most $1-2\eps$ which contains all but $O_{\mu,\eps}(n^{1-\eps})$ of the $v_1,\ldots,v_n$.  But rank has to be an integer, thus $Q$ has rank zero and is therefore just $\{0\}$.  Thus at least one of the $v_i$ is zero, a contradiction.

Now we prove (ii).  Suppose for contradiction that $\P_\mu(\bv) \geq n^{-3/2+\eps}$.  Applying Theorem \ref{theorem:strong} with $A := 3/2-\eps$ (and $\eps$ replaced by a smaller quantity $\eps'$ depending only on $\eps$) we see that there exists a symmetric GAP $Q$ of rank $d$ at most $3-2\eps$ and volume $O_{\mu,\eps}( n^{3/2 - d/2 - \eps'} )$ which contains all but $O_{\mu,\eps}(n^{1-\eps'})$ of the $v_1,\ldots,v_n$.  Since the $v_i$ are all distinct, this forces $|Q| \gg n$, which forces $d=0$, which forces more than one of the $v_i$ to be zero, a contradiction.
\end{proof}

In a similar spirit, we obtain the following variant of Theorem \ref{theorem:halasz}, which essentially asserts that equality in Theorem \ref{theorem:halasz} is only attained when the $v_i$ lie in a symmetric arithmetic progression (i.e. a symmetric rank $1$ GAP):

\begin{proposition} \label{theorem:newhalasz}  Let $n,\bv,G,\mu,l,R_l$ be as in Theorem \ref{theorem:halasz}, and let $0 < \delta, \eps < 1/2$.  Then
one of the following statements hold:
\begin{itemize}
\item $\P_{\mu} (\bv) \ll_{l, \mu, \eps, \delta} n^{-2l-1/2-\delta} R_l$.
\item All but at most $n^{1-\eps}$ of the $v_i$ lie in an symmetric arithmetic
progression of length at most $n^{2l+\delta +\eps} R_l^{-1}$.
\end{itemize}
\end{proposition}

Note that by combining this proposition with Proposition \ref{filo} and taking $\delta=\eps$ we obtain Theorem \ref{theorem:halasz} up to epsilon losses.

\begin{proof}  By shrinking $\eps$ if necessary, we may assume $\eps$ is small depending on $l, \delta$.  We may assume that $n$ is large depending on $l,\mu,\eps,\delta$, since the claim is trivial otherwise.  Finally, we may assume that
$$ \P_{\mu} (\bv) \geq n^{-2l-1/2-\delta} R_l$$
since we are clearly done otherwise.

Applying Theorem \ref{theorem:stronggeneral} with $k := n^{1/2-\eps}$ and $d = O_l(1)$ we obtain a proper symmetric GAP $Q$ of rank $r = O_l(1)$ and volume
$$ \operatorname{vol}(Q) \ll_{\mu,l,\delta,\eps} n^{2l + 1/2 + \delta + \eps} R_l^{-1} $$
such that $Q_{1/k}$ contains all but at most $n^{1-\eps}$ of the $v_i$.   Arguing as in the proof of Theorem \ref{theorem:strong}, we may assume that all dimensions of $Q$ are at least $k$.  

If $r \leq 1$ then we are done,  as $Q_{1/k}$ is an arithmetic progression having the right length (one can adjust the constant $\eps$). Now 
 assume for contradiction that $r \geq 2$.  Then (if $\eps$ is small enough) we conclude
$$ |Q_{1/k}| \ll_l k^{-2} \operatorname{vol}(Q) \ll_l n^{2l - \eps} R_l^{-1}.$$
By relabeling we may assume that the $v_1,\ldots,v_{\lfloor n/2 \rfloor}$ (say) lie in $Q_{1/k}$.  Consider the $\Theta_l(n^l)$ sums formed by taking $l$ of these $v_1,\ldots,v_{\lfloor n/2 \rfloor}$; these lie in $Q_{l/k}$, which has cardinality $O_l( n^{2l-\eps} R_l^{-1} )$.  Applying the Cauchy-Schwarz inequality, we conclude that the number of solutions to
$$ v_{i_1} + \ldots + v_{i_l} = v_{i_{l+1}} + \ldots + v_{i_{2l}}$$
with $i_1,\ldots,i_{2l} \in \{1,\ldots,\lfloor n/2\rfloor\}$, is $\gg_l n^{2l} / (n^{2l-\eps} R_l^{-1}) = n^{-\eps} R_l$.  On the other hand, this number is clearly bounded above by $R_l$, giving the required contradiction.
\end{proof}

The rest of the paper is organized as follows. In the next two sections, we recall and prove several lemmas. The proof of Theorem \ref{theorem:stronggeneral} will be presented in the last two sections of the paper. 

\section{Generalized arithmetic progressions}\label{gap-appendix}

In this section we recall the concept of a generalized arithmetic progression (GAP) and their basic properties.  
A detailed treatment of this topic can be found in \cite[Chapter 3]{TVbook}.  We will restrict our attention to symmetric GAPs.

\begin{definition}[GAPs]  Let $G$ be an additive group.  A \emph{symmetric generalized arithmetic progression} in $G$, or \emph{symmetric GAP} for short, is a quadruplet $\Q = (Q, N, v, d)$, where the \emph{rank} $\rank(\Q) = d$ is a non-negative integer, the \emph{dimensions} $N = (N_1,\ldots,N_d)$ are a $d$-tuple of positive reals, the \emph{steps} $v = (v_1,\ldots,v_d)$ are a $d$-tuple of elements of $G$, and $Q \subset G$ is the set
$$ Q = \{ \sum_{i=1}^d n_i v_i: n_i \in [-N_i,N_i] \forall i =1,\ldots,d\},$$
where $[a,b]$ denotes the set of integers between $a$ and $b$ inclusive.
We shall often abuse notation and write $Q$ for $\Q$.  For any $t > 0$, we define the \emph{dilate} $\Q_t$ of $\Q$ to be the GAP $\Q_t := (Q_t, t N, v, d)$ formed by dilating all the dimensions by $t$.   We say that $\Q$ is \emph{proper} if all the elements $n_1 v_1 + \ldots + n_d v_d$ for $n_i \in [-N_i,N_i]$ are distinct.   We say that $Q$ is \emph{$t$-proper} if $tQ$ is proper.

We define the \emph{volume} of $\Q$ to be $\operatorname{vol}(\Q) := \prod_{i=1}^d (2 \lfloor N_i \rfloor + 1)$.  Note that $|Q| \leq \operatorname{vol}(\Q)$, with equality if and only if $\Q$ is proper.
\end{definition}

If $Q$ is a GAP of rank $d$, a simple covering argument (see \cite[Lemma 3.10]{TVbook}) shows the doubling bounds
\begin{equation}\label{double}
|Q_t| \ll_d (1+t)^d |Q|
\end{equation}
for all $t > 0$.

\begin{proof}[Proof of Proposition \ref{filo}]  Let $w_1,\ldots,w_d$ be the steps of $Q$, let $N_1,\ldots,N_d$ be the dimensions, and let $\phi: \Z^d \to G$ be the homomorphism $\phi(a_1,\ldots,a_d) := a_1 w_1 + \ldots + a_d w_d$.  By hypothesis, we can write $v_i = \sum_{j=1}^d c_{ij} w_j$ for some integers $-N_j \leq c_{ij} \leq N_j$.  Then we have $S_\mu(\bv) = \phi( x )$, where $x = (x_1,\ldots,x_d) \in \Z^d$ is the random variable whose coefficients are given by
$$ x_j := \sum_{i=1}^n \eta_i c_{ij}.$$
A simple computation shows that each $x_j$ has mean zero and variance $O( N_j^2 \mu n )$, and so
$$ \E \sum_{j=1}^d |x_j|^2 / N_j^2 \ll_d \mu n.$$
By Markov's inequality, we thus conclude that
$$ \sum_{j=1}^d |x_j|^2 / N_j^2 \ll_d \mu n$$
with probability at least $1/2$ (say).  This implies that $S_\mu(\bv) \in Q_{O(\sqrt{\mu n})}$ with probability at least $1/2$, and so by the pigeonhole principle
$$ |\P_\mu(\bv)| \gg 1 / |Q_{O(\sqrt{\mu n})}|$$
and the claim follows from \eqref{double}.
\end{proof}

One can easily pass from GAPs to proper GAPs by the following lemma:

\begin{lemma} [Embedding Lemma] \cite{TVjohn} \label{lemma:embedding}
Let $Q$ be a symmetric GAP of rank $d$ in a torsion-free additive group $G$, and let $t$ be a positive constant. Then
there is a $t$-proper symmetric GAP $Q'$ of rank at most $d$ such that $Q
\subset Q' \subset Q_{O_{d,t}(1)}$ and $|Q'| \ll_{d,t} |Q|$.   If $Q$ was not already $t$-proper, one can take $Q'$ to have rank at most $d-1$.
\end{lemma}

\begin{proof} See \cite[Theorem 1.11]{TVjohn}.
\end{proof}

Recall from the homomorphism theorems that if $H, K$ are two finite subgroups of an abelian group $G$, then $|H| |K| = |H+K| |H \cap K|$.  We now establish the analogous conclusion for GAPs (cf. \cite[Exercise 2.4.7]{TVbook}):

\begin{lemma}[Intersection lemma]\label{lemma:sandwich}
Let $P$ and $Q$ be symmetric GAPs in an additive group $G$ of
rank at most $d$, then
\begin{equation}\label{pqp}
 |P\cap Q| |Q+P| = \Theta_d( |P| |Q| ).
\end{equation}
Here of course $Q+P := \{ q+p: q \in Q, p \in P \}$ denotes the sumset of $Q$ and $P$.
\end{lemma}

\begin{proof}  We recall the Ruzsa triangle inequality
$$|A-C||B| \le |A-B| |B-C|$$
for finite non-empty sets $A,B,C \subset G$ (see e.g. \cite[Lemma 2.6]{TVbook}); this follows from the fact that any element $a-c $ with $a \in A$
and $c \in C$ has at least $b$ representations of the form $a-c = (a-b)+ (b-c) $ with $b \in B$.  Applying this with $A=P, C=Q, B=P \cap Q$ we obtain
$$|P-Q||P \cap Q| \le | P- (P\cap Q)| | (P\cap Q) -Q | \leq |2P| |2Q|$$
where we use the symmetry of $P, Q$.  But from \eqref{double} we have $|2P| \ll_d |P|$, $|2Q| \ll_d |Q|$, which gives the upper bound in \eqref{pqp}.

Now we turn to the lower bound. By reducing $d$ if necessary, we can assume that the
dimensions of both $P$ and $Q$ are divisible by two, thus $P = P_{1/2} - P_{1/2}$ and $Q = Q_{1/2} - Q_{1/2}$.  Now we recall the inequality
$$ |A| |B| \le |(A-A) \cap (B-B)| |A+B|$$
for finite non-empty $A,B \subset G$ (cf. \cite[Corollary 2.10]{TVbook}), which follows by combining the identity
$$|A||B| = |\{(a,b)|a \in A, b \in B \}|=  \sum_{x \in A+B} |\{(a,b)| a \in A, b \in B, a+b = x\} | $$
with the inequality
$$ |\{(a,b)| a \in A, b \in B, a+b = x\} | \leq |(A-A) \cap (B-B)|$$
for all $x \in G$ (which follows from the observation that if $a,a' \in A$ and $b,b' \in B$ are such that $a+b=a'+b'=x$, then $a-a'=b-b'$ lies in $(A-A) \cap (B-B)$).  Applying this inequality with $A=P_{1/2}$ and $B=Q_{1/2}$ and using \eqref{double}, one obtains the claim.
\end{proof}

\section{Arithmetic on words} \label{fourier}

In this section, we recall some tools developed  earlier in
\cite{TVcond}, which were used to prove Theorem \ref{theorem:weak}
and will be useful here as well.

For our purpose, it  is convenient to think of $\bv = (v_1, \dots, v_n)$ as a word,
obtained by concatenating the $v_i$:
$$\bv = v_1 v_2 \dots v_n. $$
This will allow us to perform several operations such as
concatenating, truncating and repeating. For instance, if $\bv
=v_1\dots v_n$ and $\bw= w_1\dots w_m$, then
$$\P_{\mu} (\bv \bw) = \max _{a \in Z} \Big(\sum_{i=1}^n \eta_i^{\mu} v_i
+ \sum_{j=1}^m \eta_{n+j} ^{\mu} w_j=a \Big)$$ \noindent where
$\eta_{k}^{\mu}, 1\le k \le n+m$ are  i.i.d copies of $\eta^{\mu}$.
Furthermore, we use $\bv^{[k]}$ to denote the concatenation of $k$
copies of $\bv$.  

We will need to generalize the concentration probabilities $\P_\mu(\bv)$ as follows. For finite non-empty set $Q \subset G$, define
\begin{equation}\label{pbq-def}
\P_\mu(\bv;Q) := \sup_{a \in G} \P( S^\mu(\bv) = a + q - q' )
\end{equation}
where $q, q'$ are independently chosen uniformly at random from $Q$.  Note that $\P_\mu(\bv;Q) = \P_\mu(\bv)$ if $Q$ is a singleton set.

Since $\P( a+q-q'=x) \leq 1/|Q|$ for any fixed $a,q', x$, a simple conditioning argument reveals the crude bound
\begin{equation}\label{pbound}
 \P_\mu(\bv;Q)  \leq 1/|Q|.
\end{equation}

We have the following basic properties of the $\P_\mu(\bv)$ and $\P_\mu(\bv;Q)$:

\begin{lemma}\label{lemma:properties} Let $\bv = v_1\ldots v_n$ be a word $v_1,\ldots,v_n$ in a torsion-free additive group $G$, and let $Q \subset G$ be a finite non-empty set.  Then the following properties
hold.
\begin{itemize}

\item[(i)] $\P_\mu(\bv;Q)$ is  invariant under permutations of $\bv$.

\item[(ii)] For any words  $\bv, \bw$

\begin{equation} \notag \P_\mu(\bv \bw;Q) \leq \P_\mu(\bv;Q).
\end{equation}

\item[(iii)] For any $0 < \mu \le 1$, any $0 < \mu' \leq \mu/4$, and any word $\bv$,
\begin{equation} \notag \label{14} \P_{\mu}(\bv;Q) \leq \P_{\mu'}(\bv;Q).
\end{equation}

\item[(iv)]  For any number $0 < \mu \leq 1/2$ and any word $\bv$,
\begin{equation} \notag \label{pmub-2} \P_\mu(\bv;Q) \leq \P_{\mu/k}(\bv^{[k]};Q).
\end{equation}

\item[(v)] For any number $0 < \mu \leq 1/2$ and any words $\bv, \bw_1,
\dots, \bw_m$ we have
\begin{equation} \notag \label{pmub-3}
\P_\mu(\bv \bw_1 \ldots \bw_m;Q) \leq \left( \prod_{j=1}^m \P_\mu(\bv
\bw_j^{[m]};Q) \right)^{1/m}.
\end{equation}

\item[(vi)]  For any number $0 < \mu \leq 1/2$ and any words $\bv,
\bw_1, \dots, \bw_m$, there is an index $1\le j \le m$ such that
\begin{equation}\notag \label{pmub-3a}
\P_\mu(\bv \bw_1 \ldots \bw_m;Q) \leq \P_\mu(\bv \bw_j^{[m]};Q).
\end{equation}
\end{itemize}
\end{lemma}

\begin{proof}  When $G=\Z$ and $Q$ is a singleton, this is \cite[Lemma 5.1]{TVcond}.  When $G=\Z$ and $Q$ is not a singleton, the claim can be established by repeating the proof of \cite[Lemma 5.1]{TVcond}, using the Fourier identity
$$ \P( S^\mu(\bv) = a + q - q' ) = \int_{0}^{1 } e(-a t) |\E (e(qt))|^{2 } \prod_{i=1}^{n} (1-\mu + \mu \cos 2\pi v_{i }t ) \,\, dt$$
in place of
$$ \P( S^\mu(\bv) = a ) = \int_{0}^{1 } e(-a t) \prod_{i=1}^{n} (1-\mu + \mu \cos 2\pi v_{i }t ) \,\, dt;$$
we omit the details.  
(Here and later, $e(x)$ denotes $ \exp(2\pi \sqrt{-1} x) $.) Finally, the generalization to arbitrary torsion-free $G$ can be accomplished by using Freiman isomorphisms (see \cite[Lemma 5.25]{TVbook}).
\end{proof}

Note that for fixed $0 < \mu < 1$, a random walk $S^\mu( v^{[k^2]} )$ is roughly uniformly distributed on the progression $[-k,k] v := \{ j v: j \in \Z, -k \leq j \leq k \}$, thanks to the central limit theorem.  (Here $v^{[k^2]}$ is the word $v$ repeated $k^2$ times.) The following lemma can be viewed as a formalization of this intuition.

\begin{proposition}[Comparison of random walks]\label{cor:prob3}  Let $0 < \mu \leq 1/2$, let $\bv = (v_1,\ldots,v_n)$ be a tuple in a torsion-free additive group $G$, let $v_0 \in G$, and let $k \geq 1$.  Let $Q$ be a symmetric GAP in $\Z$ of rank $d$.  Then
$$ \P_\mu( \bv v_0^{[k^2]}; Q ) \ll_{\mu,d} \P_\mu( \bv; Q + [-k,k]v_0 ).$$
\end{proposition}

\begin{proof}  Fix $\mu,d$; we allow all implied constants to depend on these quantities.  
By definition,
$$ \P_\mu( \bv v_0^{[k^2]}; Q ) = \P( S^\mu(\bv) + X v_0 = q-q' )$$
where $q,q'$ are independent random variables uniformly distributed in $Q$, and $X := \sum_{i=1}^{k^2} \xi_i^\mu$.  A direct computation using Stirling's formula shows that
$$ \P( X = m ) \ll k^{-1} \exp( - \Omega(|m|/k) )$$
for all $m \in \Z$, thus
$$ \P_\mu( \bv v_0^{[k^2]}; Q ) \ll k^{-1} \sum_{m \in \Z} \exp( - \Omega(|m|/k) ) \P( S^\mu(\bv) + m v_0 = q-q' ).$$
This implies that
$$ \P_\mu( \bv v_0^{[k^2]}; Q ) \ll k^{-1} \sum_{m \in \Z} \exp( - \Omega(|m|/k) ) \P( S^\mu(\bv) + m v_0 = q-q'+jv_0-j'v_0 )$$
where $j, j'$ are drawn uniformly at random from $[-k,k]$, independently of each other and of $q, q'$.  
It therefore suffices to show that
$$ \P( S^\mu(\bv) = a+q-q'+jv_0-j'v_0 ) \ll \P_\mu( \bv; Q + [-k,k]v_0 )$$
for all $a \in G$.

The random variable $q+jv_0$ is supported in $Q + [-k,k]v_0$.  If it were distributed uniformly in this set, we would be done.  It is not quite uniform, nevertheless we can compare it to the uniform distribution as follows.  Given any $x \in Q + [-k,k]v_0$, we have
$$ \P( q+jv_0 = x ) = \frac{1}{|Q| |[-k,k]v_0|} |Q \cap (x - [-k,k]v_0)|.$$
Since $|A| \leq |A-A|$, we see that
$$ |Q \cap (x - [-k,k]v_0)| \leq |(Q-Q) \cap ([-2k,2k]v_0)|$$
and so by Lemma \ref{lemma:sandwich} and \eqref{double}
$$ |Q \cap (x - [-k,k]v_0)| \ll \frac{|Q| |[-k,k]v_0|}{|Q + [-k,k]v_0|}$$
and thus
$$ \P( q+jv_0 = x ) \ll \frac{1}{|Q + [-k,k]v_0|}.$$
Thus the probability distribution of $q+jv_0$ is majorized by a constant multiple of the uniform distribution on $Q+[-k,k]v_0$, and the claim follows.
\end{proof}


\section{The algorithm} \label{algorithm}

We begin the proof of Theorem \ref{theorem:stronggeneral}.  By Lemma \ref{lemma:properties} we may assume that $\mu \leq 1/2$.  
Fix $d, \eps, \mu, n, k, \bv, G$ as in that theorem; we assume that $n, k$ are sufficiently large depending on $d,\eps,\mu$.  We let $K \geq 1$ be a large number depending on $d,\eps,\mu$, and then let $C_0$ be an even larger number depending on $d,\eps,\mu,K$.  We assume that \eqref{pi2} holds.  

In this section, we describe an algorithm which takes $\bv$ as  input and produces, as output, a symmetric GAP $Q$ as claimed by Theorem \ref{theorem:stronggeneral}.  A key concept is that of a \emph{bad} element with respect to a symmetric GAP.

\begin{definition}[Bad element]  Let $K \geq 1$, $x \in G$, and let $Q$ be a symmetric GAP in $G$.  We say that $x$ is \emph{bad} with respect to a symmetric GAP $Q$ if
$$|Q  + [-k,k] x | \ge K |Q|$$
and \emph{good} otherwise.
\end{definition}

We will also need the generalized concentration probabilities $\P_\mu(\bv;Q)$ defined in \eqref{pbq-def}.  We now consider the following algorithm that generates words $\bv^i$ and symmetric GAPs $Q_i$ for various $i=0,1,2,\ldots$:

\vskip2mm

{\bf Step 0.} Set  $\bv^{0}= \bv, Q_{0}:= \{0\}$.

{\bf Step $i+1$.} Count the number of elements of $\bv^i$  which are
bad with respect to $Q_i$.

{\it Case 1.}  If this number is less than $k^2$  then STOP.

{\it Case 2.}  If this number is at least $k^2$, we can assume
(without loss of generality) that the last $k^{2}$ coordinates of
$\bv^{i}$ are bad. Let $\bv^{i+1}$ be the vector obtained from
$\bv^{i}$ by truncating these bad coordinates. By Lemma
\ref{lemma:properties}(vi), there is some value $v_0$ among the bad coordinates such
that
$$P_{\mu} (\bv^{i+1} v_0 ^{[k^{2}]}; Q_{i}) \ge  \P_{\mu}(\bv^{i}; Q_{i}). $$
Set $r_i := \operatorname{rank}(Q_i)$ and $Q_{i+1}' : = Q_{i} + [-k,k] v_{0}$, thus $Q'_{i+1}$ is a GAP with rank $r_i+1$. If $Q'_{i+1}$ is proper,
then set $Q_{i+1}:= Q'_{i+1}$. If it is not proper, then use Lemma
\ref{lemma:embedding} to embed it into a proper symmetric GAP of
rank at most $r_i$ and volume $O_{r_i}(|Q_{i+1}'|)$;. Call this
proper GAP $Q_{i+1}$.   CONTINUE  to Step $i+2$.

\vskip2mm

Notice that by the algorithm the $Q^{i}$ are symmetric GAPs at every step.

\section{Analysis of the algorithm} \label{analysis1}

For each $i$ that occurs in the algorithm, we define the rank
$$ r_i := \operatorname{rank}(Q_i)$$
and the potential
$$ F_i := |Q_i| \P_\mu(\bv^i; Q_i).$$
Initially we have
\begin{equation}\label{r0f0}
 r_0 = 0; \quad F_0 = \P_\mu(\bv) \geq C_0 k^{-d}; \quad |Q_0| = 1.
\end{equation}

We now record how $r_i$, $F_i$, and $Q_i$ evolve with the algorithm.  
We say that Step $i+1$  is \emph{proper} if $Q'_{i+1}$ is proper. 

\begin{lemma}\label{grow}  Let Step $i+1$ be a step that occurs in the algorithm.
\begin{itemize}
\item[(i)] We have $r_{i+1}=r_i+1$ if Step $i+1$ is proper, and $r_{i+1} \leq r_i$ otherwise.
\item[(ii)] We have $|Q_{i+1}| \geq K |Q_i|$.  If Step $i+1$ is proper, we can improve this to $|Q_{i+1}| \geq k |Q_i|$.
\item[(iii)] We have $F_{i+1} \gg_{r_i,\mu} K F_i$.  If Step $i+1$ is proper, we can improve this to $F_{i+1} \gg_{r_i,\mu} k F_i$.
\end{itemize}
\end{lemma}

\begin{proof} The first two claims are clear from construction.  To prove the third claim, we observe from Proposition \ref{cor:prob3} that
\begin{equation}\label{pmuq}
\P_{\mu} (\bv^{i+1}; Q_{i+1} ) \gg_{r_i,\mu} \P_{\mu} (\bv^{i} ; Q_i);
\end{equation}
the claim (iii) now follows from (ii).
\end{proof}

\begin{corollary}\label{term} The algorithm has at most $d-1$ proper steps, and terminates in $O( d \log_K k )$ steps.
\end{corollary}

\begin{proof}  Suppose for contradiction that there were at least $d$ proper steps.  Let $1 \leq i_1 < \ldots < i_{d}$ be the first $d$ proper steps.  By Lemma \ref{grow}(i) and \eqref{r0f0}, the ranks $r_i$ are bounded by $d$ for all $i \leq i_d$.  From Lemma \ref{grow}(iii), we have $F_{i_d} \geq_{d,\mu} k^d F_0$ if $K$ is large enough; on the other hand, from \eqref{pbound} we have $F_{i_d} \leq 1$.  This contradicts \eqref{r0f0} if $C_0$ is large enough.

Now that there are at most $d-1$ proper steps, $r_i \leq d-1$ for all $i$.  By Lemma \ref{grow}(iii), we thus have $F_{i+1} \geq \sqrt{K} F_i$ for all $i$ if $K$ is large enough.  On the other hand, from \eqref{pbound} we have $F_i \leq 1$ for all $i$.  Applying \eqref{r0f0}, we conclude that the algorithm terminates in $O(d \log_K k)$ steps as claimed.
\end{proof}

Let $\bv^{T}$ and $Q_{T}$ be the vector and GAP at the stopping
time $T = O(d \log_K k)$.  

\begin{lemma} $Q_T$ has rank at most $d-1$ and
$$ |Q_T| \leq k^{\eps/2} \P_\mu(\bv)^{-1}$$
\end{lemma}

\begin{proof} The rank bound follows from \eqref{r0f0}, Lemma \ref{grow}(i) and Corollary \ref{term}. 

As proved above, $Q_T$ has rank at most $d-1$. We next prove
that it has small cardinality.  Iterating \eqref{pmuq} starting from \eqref{r0f0}, we see that
$$ \P_{\mu} (\bv^{T}; Q_{T}) \geq \Omega_{d,\mu}(1)^T \P_\mu(\bv).$$
Combining this with \eqref{pbound} and the bound $T = O(d \log_K k)$ we conclude
$$ |Q_T| \leq \exp( O_{d,\mu}(\log_K k) ) \P_\mu(\bv)^{-1}$$
and the claim follows by taking $K$ sufficiently large.
\end{proof}

By construction, all but $O(d k^2 \log_K k) = O( k^2 \log k )$ of the $v_1,\ldots,v_n$ are good relative to $Q_T$.  To exploit this we use
 
\begin{lemma} \label{lemma:good}  Suppose that $x \in G$ is good relative to a symmetric GAP $Q$ of rank $r$.  Then there exists a proper symmetric GAP $Q'$ of rank at most $r$ containing $Q$ and volume $|Q'| \ll_{K,r} |Q|$ such that $Cx \in Q'_{C/k}$, where $C \geq 1$ is an integer depending only on $K$ and $r$.
\end{lemma}

\begin{proof}  The $|Q| |[-k,k]|$ sums $q+jx$ with $q \in Q$ and $j \in [-k,k]$ lie in the set $|Q + [-k,k]x|$, which has cardinality at most $K|Q|$ by hypothesis.  By Cauchy-Schwarz, we conclude that there are $\gg_K |Q| k^2$ quadruplets $(q,q',j,j') \in Q \times Q \times [-k,k] \times [-k,k]$ such that $q+jx=q'+j'x$.  By the pigeonhole principle, we conclude that the set $A := \{ j \in [-2k,2k]: jx \in Q-Q \}$ has cardinality $|A| \gg_K k$.   Applying a result of S\'ark\H{o}zy \cite{sar} (see also \cite{lev}, \cite{sv}, or \cite[Chapter 12]{TVbook}) we conclude that there exists a positive integer $K_1 = O_K(1)$ such that the iterated sumset $K_1 A$ contains an arithmetic progression of length $l = \Theta_K(k)$ and positive integer step $a = O_K(1)$.  We conclude that $[-l,l]ax \in Q_{4K_1}$.

At present, $l$, $a$, and $K_1$ are all dependent on $x$.  But $K_1, a$ are bounded by $O_K(1)$, and $l$ is bounded from below by $\Omega_K(k)$.  Thus, by taking the gcd over all possible values of $a$, one may assume that $l, a, K_1$ are independent of $x$.

By Lemma \ref{lemma:embedding}, we can place $Q_{4K_1}$ inside a $2$-proper symmetric GAP $Q'$ of rank at most $r$ and volume $O_{K,r}(|Q|)$, thus $[-l,l]ax \in Q'$.  Write $N'_1,\ldots,N'_{r'}$ for the dimensions of $Q'$.  Since $Q'$ is $2$-proper, the obvious map $\phi: [-N'_1,N'_1] \times \ldots \times [-N'_r,N'_r] \to Q'$ is a Freiman isomorphism of order $2$ (see \cite[Section 5.3]{TVbook}), and $\phi^{-1}([-l,l]ax) = [-l,l]\phi^{-1}(ax)$ is also an arithmetic progression.  From this we see that
$$ \phi^{-1}(ax) \in [-N'_1/l,N'_1/l] \times \ldots \times [-N'_r/l,N'_r/l]$$
and thus
$$ ax \in Q'_{1/l}$$
and the claim follows.
\end{proof}

Through the proof of  this lemma, we see that all but at most $O(k^2 \log k)$ of the $v_i$ are such that $Cv_i \in (Q_T)_{C/k}$.  By Lemma \ref{lemma:embedding}, we may place $(Q_T)_C$ inside a proper symmetric GAP $Q$ of rank at most $d-1$ and volume
$$|Q| \le k^{\eps} \P_{\mu} (\bv)^{-1}. $$
Since $(Q_T)_C \subset Q'$, we have $(Q_T)_{C/k} \subset Q_{1/k}$, and Theorem \ref{theorem:stronggeneral} follows.

{\it Acknowledgement.} We would like to thank the referees for their  careful reading and useful remarks.

\end{document}